\documentclass[12pt]{article}
\usepackage{amsmath,amssymb,amsthm}
\usepackage{hyperref}
\usepackage{datetime}
\usepackage{units}
\usepackage{color}
\usepackage[T1]{fontenc}
\usepackage[utf8]{inputenc}
\usepackage{authblk}
\usepackage{bm,latexsym,mathrsfs,enumerate}
\usepackage{color}

\title{ Infinite arctangent sums involving Fibonacci and Lucas numbers\thanks{%
MSC 2010: 11B39, 11Y60}}
\author[]{Kunle Adegoke\thanks{adegoke00@gmail.com }}
\affil{Department of Physics and Engineering Physics, \mbox{Obafemi Awolowo University, Ile-Ife, 220005 Nigeria}}

\theoremstyle{plain}
\numberwithin{equation}{section}
\newtheorem{thm}{THEOREM}[section]

\begin{document}
\date{}
\maketitle
\begin{abstract}
\noindent Using a straightforward elementary approach, we derive numerous infinite arctangent summation formulas involving Fibonacci and Lucas numbers. While most of the results obtained are new, a couple of celebrated results appear as particular cases of the more general formulas derived here.
\end{abstract}
\tableofcontents

\section{Introduction}

It is our goal, in this work, to derive infinite arctangent summation formulas involving Fibonacci and Lucas numbers. The results obtained will be found to be of a more general nature than one finds in earlier literature. 

\bigskip

Previously known results containing arctangent identities and/or infinite summation involving Fibonacci numbers can be found in references \cite{bragg, hayashi, hoggattv, mahon85,melham95} and references therein.

\bigskip

In deriving the results in this paper, the main identities employed are the trigonometric addition formula
\begin{equation}
\tan^{-1}\left\{\frac{\lambda {(y-x)}}{xy+\lambda^2}\right\}=\tan^{-1}\frac{\lambda}{x}-\tan^{-1}\frac{\lambda}{y}\,,\label{equ.arctan}
\end{equation}
which holds for $\lambda\in\mathbb{R}$ such that either $xy>0$ or $xy<0$ and $\lambda^2<-xy$, and the following identities which resolve products of Fibonacci and Lucas numbers

\begin{subequations}\label{equ.v128kwt}
\begin{eqnarray}\label{equ.xj19lbd}
F_{u-v}F_{u+v}=F^2_u-(-1)^{(u-v)}F_v^2\,,\label{equ.u38jrgd}\\
L_{u-v}L_{u+v}=L_{2u}+(-1)^{(u-v)}L_{2v}\,,\label{equ.itu9olz}\\
L_uF_v=F_{v+u}+(-1)^uF_{v-u}\,,\label{equ.wogl40r}\\
F_uL_v=F_{v+u}-(-1)^uF_{v-u}\,,\label{equ.fen1v2f}\\
L_{u}L_{v}=L_{u+v}+(-1)^uL_{v-u}\,,\label{equ.l5tzgky}\\
5F_{u-v}F_{u+v}=L_{2u}-(-1)^{(u-v)}L_{2v}\label{equ.kt19j3t}\,.
\end{eqnarray}
\end{subequations}

Also we shall make repeated use of the following identities connecting Fibonacci and Lucas numbers:

\begin{subequations}\label{equ.u4fy3q3}
\begin{eqnarray}\label{equ.wsr7wzv}
F_{2u}=F_uL_u\,,\\
L_{2u}-2(-1)^u=5F_u^2\,,\\
5F_u^2-L_u^2=4(-1)^{(u+1)}\,,\\
L_{2u}+2(-1)^u=L_u^2\,.
\end{eqnarray}
\end{subequations}

Identities~\eqref{equ.v128kwt} and~\eqref{equ.u4fy3q3} or their variations can be found in~\cite{basin2, dunlap, howard}.

\bigskip

On notation, $G_i$, $i$ integers, denotes generalized Fibonacci numbers defined through the second order recurrence relation \mbox{$G_i=G_{i-1}+G_{i-2}$}, where two boundary terms, usually $G_0$ and $G_1$, need to be specified. When $G_0=0$ and $G_1=1$, we have the Fibonacci numbers, denoted $F_i$, while when $G_0=2$ and $G_1=1$, we have the Lucas numbers, denoted $L_i$.

\bigskip

Throughout this paper, the principal value of the arctangent function is assumed.

\bigskip

Interesting results obtained in this paper, for integers $k$, $j\ne 0$ and $p$ include
\begin{align*}
& \sum_{r = 1}^\infty  {\tan ^{ - 1} \left\{ {\frac{{F_j^2 L_j L_{4jr} }}{{F_{4jr}^2  - F_{2j}^2  + F_j^2 }}} \right\} = \tan ^{ - 1} \left( {\frac{1}{{L_j }}} \right)}\,, \sum_{r = 1}^\infty  {\tan ^{ - 1} \left\{ {\frac{{L_j^2 F_j L_{4jr} }}{{F_{4jr}^2  - F_{2j}^2  + L_j^2 }}} \right\} = \tan ^{ - 1} \left( {\frac{1}{{F_j }}} \right)}\,,\\
&\sum_{r=1}^\infty\tan^{-1}\left\{\frac{F_{2j}^2L_{4jr+2j}}{F_{4jr+2j}^2}\right\}=\tan^{-1}\left(\frac{1}{L_{2j}}\right)\,,\sum_{r = 1}^\infty  {\tan ^{ - 1} \left\{ {\frac{{F_{2(2j - 1)} }}{{F_{4jr - 2r + 2j - 1} }}} \right\}}  = \tan ^{ - 1} \left( {\frac{{1}}{{L_{2j - 1} }}} \right)\,,\\
&\\
&\sum_{r = 1}^\infty  {\tan ^{ - 1} \left\{ {\frac{{F_{2(2j - 1)} }}{{F_{4jr - 2r + 1} }}} \right\}}  = \tan ^{ - 1} \left( {\frac{{F_{2j - 1} }}{{F_{2j} }}} \right),\,\sum_{r = 1}^\infty  {\tan ^{ - 1} \left\{ {\frac{{1}}{{F_{2r + 2k - 1} }}} \right\}}  = \tan ^{ - 1} \left( {\frac{1}{{F_{2k} }}} \right)\,,
\end{align*}

\begin{align*}
&\sum_{r = p}^\infty  {\tan ^{ - 1} \left\{ {\frac{1}{5}\frac{{L_{2r} }}{{F_{2r}^2 }}} \right\}}  = \tan ^{ - 1} \left( {\frac{1}{{L_{2p - 1} }}} \right)\,,\sum_{r=1}^\infty \tan^{-1}\left\{\frac{F_{4j}}{F_{4jr-1}}\right\}=\tan^{-1}\left(\frac{L_{2j}}{L_{2j-1}}\right)\,.
\end{align*}


We also obtained the following special values
\begin{align*}
&\sum_{r=1}^\infty\tan^{-1}\left\{\frac{L_{4r-2}}{F_{4r-2}^2}\right\}=\frac{\pi}{2}\,,\quad {\sum_{r=1}^\infty\tan^{-1}\left\{\frac{L_{4r}}{F_{4r}^2}\right\}=\frac{\pi}{4}}\,,\quad {\sum_{r = 1}^\infty  {\tan ^{ - 1} \left\{ {\frac{\sqrt{35}\, L_{4r - 2} }{{L_{2(4r - 2)} }}} \right\}}  = \frac{\pi}{2}}\,,\\
&{\sum_{r = 1}^\infty  {\tan ^{ - 1} \left\{ {\frac{{\sqrt 3 L_{2r} }}{{L_{4r} }}} \right\}}  = \frac{\pi }{3}}\,,\quad {\sum_{r = 1}^\infty  {\tan ^{ - 1} \left\{ {\frac{1}{5}\frac{{L_{2r} }}{{F_{2r}^2 }}} \right\}}  = \frac{\pi }{4}}\,,\quad {\sum_{r = 1}^\infty  {\tan ^{ - 1} \left\{ {\frac{{\sqrt 5 }}{{L_{2r} }}} \right\}}  = \tan ^{ - 1} \sqrt 5}\,,\\
&{\sum_{r = 1}^\infty  {\tan ^{ - 1} \left\{ {\frac{\sqrt{35}\, L_{4r} }{{L_{8r} }}} \right\}}  = \sqrt{\frac{7}{5}}}\,,\, {\sum_{r = 1}^\infty  {\tan ^{ - 1} \left\{ {\frac{{\sqrt 5 F_{2r - 1} }}{{L_{2r - 1}^2 }}} \right\}}  = \frac{\pi}{2}}\,,\\
&{\sum_{r=1}^\infty\tan^{-1}\left(\frac{5\sqrt {7} F_{4r-1}}{L_{2(4r-1)}}\right)=\tan^{-1}\sqrt 7}\,,\quad\sum_{r=1}^\infty\tan^{-1}\left\{\frac{L_{4r+2}}{F_{4r+2}^2}\right\}=\tan^{-1}\left(\frac{1}{3}\right)\,.
\end{align*}


\section{Preliminary result}

Taking $x=G_{mr + n - m}$ and $y=G_{mr + n}$ in the arctangent addition formula, identity~\eqref{equ.arctan}, gives

\begin{equation}\label{equ.jeavdxa}
\tan ^{ - 1} \left\{ {\frac{{\lambda (G_{mr + n}  - G_{mr + n - m} )}}{{G_{mr + n} G_{mr + n - m}  + \lambda ^2 }}} \right\} = \tan ^{ - 1} \left( {\frac{\lambda }{{G_{mr + n - m} }}} \right) - \tan ^{ - 1} \left( {\frac{\lambda }{{G_{mr + n} }}} \right)\,.
\end{equation}

Summing each side of identity~\eqref{equ.jeavdxa} from $r=p\in\mathbb{Z}$ to $r=N\in\mathbb{Z^+}$ and noting that the summation of the terms on the right hand side telescopes, we obtain

\begin{equation}
\sum\limits_{r = p}^N {\tan ^{ - 1} \left\{ {\frac{{\lambda (G_{mr + n}  - G_{mr + n - m} )}}{{G_{mr + n} G_{mr + n - m}  + \lambda ^2 }}} \right\}}  = \tan ^{ - 1} \left( {\frac{\lambda }{{G_{mp + n - m} }}} \right) - \tan ^{ - 1} \left( {\frac{\lambda }{{G_{mN + n} }}} \right)\,.
\end{equation}

Now taking limit as $N\to\infty$, we have

\subsection*{Lemma.}
For $\lambda\in\mathbb{R}$, $n,m,p\in\mathbb{Z}$, $m\ne 0$ holds

\begin{equation}\label{equ.ahydik3}
\sum_{r = p}^\infty  {\tan ^{ - 1} \left\{ {\frac{{\lambda  {(G_{mr + n} - G_{mr + n - m} }) }}{{G_{mr + n} G_{mr + n - m}  + \lambda ^2 }}} \right\}}=\tan ^{ - 1} \left( {\frac{\lambda }{{G_{mp+n-m} }}} \right)\,.
\end{equation}

\section{Main Results}

\subsection{$G\equiv F$ in identity~\eqref{equ.ahydik3}, that is, $G_0=0$, $G_1=1$}

Choosing $m=4j$ and $n=2k+2j$ and using identities~\eqref{equ.u38jrgd} and \eqref{equ.fen1v2f} we prove
\begin{thm}\label{thm.xrc79et}

For $\lambda\in\mathbb{R}$, $j,k,p\in\mathbb{Z}$ and $j\ne 0$ holds
\begin{equation}\label{equ.o14ht74}
\sum_{r=p}^\infty\tan^{-1}\left\{\frac{\lambda F_{2j}L_{4jr+2k}}{F_{4jr+2k}^2-F_{2j}^2+\lambda^2}\right\}=\tan^{-1}\left(\frac{\lambda}{F_{4jp+2k-2j}}\right)\,,
\end{equation}
\end{thm}

while taking $m=4j-2$ and $n=2k+2j-2$ and using identities~\eqref{equ.u38jrgd} and \eqref{equ.wogl40r} we prove

\begin{thm}\label{thm.kc5pf5v}

For $\lambda\in\mathbb{R}$ and $j,k,p\in\mathbb{Z}$ holds

\begin{equation}\label{equ.wz6niw6}
\sum_{r = p}^\infty  {\tan ^{ - 1} \left\{ {\frac{{\lambda L_{2j - 1} F_{4jr - 2r + 2k - 1} }}{{F_{4jr - 2r + 2k - 1}^2  - F_{2j - 1}^2  + \lambda ^2 }}} \right\}}  = \tan ^{ - 1} \left( {\frac{\lambda }{{F_{4jp - 2p + 2k - 2j} }}} \right)\,.
\end{equation}

\end{thm}

\subsection{$G\equiv L$ in identity~\eqref{equ.ahydik3}, that is, $G_0=2$, $G_1=1$}

Choosing $m=4j$ and $n=2k+2j-1$ and using identities~\eqref{equ.itu9olz} and \eqref{equ.kt19j3t} we prove

\begin{thm}\label{thm.xi2xcnh}

For $\lambda\in\mathbb{R}$, $j,k,p\in\mathbb{Z}$ and $j\ne 0$ holds

\begin{equation}\label{equ.ywa2m9d}
\sum_{r=p}^\infty\tan^{-1}\left(\frac{5\lambda F_{2j}F_{4jr+2k-1}}{L_{8jr+4k-2}-L_{4j}+\lambda^2}\right)=\tan^{-1}\left(\frac{\lambda}{L_{4jp+2k-2j-1}}\right)\,,
\end{equation}

\end{thm}

while taking $m=4j-2$ and $n=2k+2j-1$ and using identities~\eqref{equ.itu9olz} and \eqref{equ.l5tzgky} we prove

\begin{thm}\label{thm.vqhpot1}

For $\lambda\in\mathbb{R}$ and $j,k,p\in\mathbb{Z}$ holds

\begin{equation}\label{equ.ocyrpi4}
\sum_{r=p}^\infty\tan^{-1}\left(\frac{\lambda L_{2j-1}L_{4jr-2r+2k}}{L_{8jr-4r+4k}-L_{4j-2}+\lambda^2}\right)=\tan^{-1}\left(\frac{\lambda}{L_{4jp-2p+2k-2j+1}}\right)\,.
\end{equation}

\end{thm}

\section{Corollaries and special values}

Different combinations of the parameters $\lambda$, $j$, $k$ and $p$ in the above theorems yield a variety of interesting particular cases. In this section we will consider some of the possible choices.

\subsection{Results from Theorem~\ref{thm.xrc79et}}






\subsubsection{$\lambda=F_{j}$, $p=1$ and $k=0$ in identity~\eqref{equ.o14ht74}}

The choice $\lambda=F_{j}$, $p=1$ and $k=0$ in identity~\eqref{equ.o14ht74} gives

\begin{equation}
\sum_{r = 1}^\infty  {\tan ^{ - 1} \left\{ {\frac{{F_j^2 L_j L_{4jr} }}{{F_{4jr}^2  - F_{2j}^2  + F_j^2 }}} \right\} = \tan ^{ - 1} \left( {\frac{1}{{L_j }}} \right)}\,. 
\end{equation}

Thus, at $j=1$, we obtain the special value

\begin{equation}\label{equ.s6rsk7b}
\boxed{\sum_{r=1}^\infty\tan^{-1}\left\{\frac{L_{4r}}{F_{4r}^2}\right\}=\frac{\pi}{4}}\,.
\end{equation}

\subsubsection{$\lambda=L_{j}$, $p=1$ and $k=0$ in identity~\eqref{equ.o14ht74}}

The above choice gives

\begin{equation}\label{equ.w0riexv}
\sum_{r = 1}^\infty  {\tan ^{ - 1} \left\{ {\frac{{L_j^2 F_j L_{4jr} }}{{F_{4jr}^2  - F_{2j}^2  + L_j^2 }}} \right\} = \tan ^{ - 1} \left( {\frac{1}{{F_j }}} \right)}\,. 
\end{equation}

At $j=1$, identity\eqref{equ.s6rsk7b} is reproduced, while at $j=2$ we have the special value

\begin{equation}\label{equ.zxs4hwa}
\boxed{\sum_{r=1}^\infty\tan^{-1}\left\{\frac{9L_{8r}}{F_{8r}^2}\right\}=\frac{\pi}{4}}\,.
\end{equation}

Note that identities~\eqref{equ.s6rsk7b} and \eqref{equ.zxs4hwa} are special cases of identity~\eqref{equ.t70h9q0} below, at $j=1$ and $j=2$, respectively.

\subsubsection{$\lambda=F_{2j}$, $k=j$ and $p=0$ in identity~\eqref{equ.o14ht74}}

This choice gives 
\begin{equation}\label{equ.zr8gn0u}
\sum_{r=1}^\infty\tan^{-1}\left\{\frac{F_{2j}^2L_{4jr-2j}}{F_{4jr-2j}^2}\right\}=\frac{\pi}{2}\,,
\end{equation}

which, at $j=1$, gives the special value

\begin{equation}\label{equ.bmwsq09}
\boxed{\sum_{r=1}^\infty\tan^{-1}\left\{\frac{L_{4r-2}}{F_{4r-2}^2}\right\}=\frac{\pi}{2}}\,.
\end{equation}

\subsubsection{$\lambda=F_{2j}$ and $p=1$ in identity~\eqref{equ.o14ht74}}

This choice gives

\begin{equation}\label{equ.bxwz09x}
\sum_{r=1}^\infty\tan^{-1}\left\{\frac{F_{2j}^2L_{4jr+2k}}{F_{4jr+2k}^2}\right\}=\tan^{-1}\left(\frac{F_{2j}}{F_{2j+2k}}\right)\,.
\end{equation}

At $k=0$ in identity~\eqref{equ.bxwz09x} we have

\begin{equation}\label{equ.t70h9q0}
\sum_{r = 1}^\infty  {\tan ^{ - 1} \left\{ {\frac{{F_{2j}^2 L_{4jr} }}{{F_{4jr}^2 }}} \right\}}  = \frac{\pi }{4}\,.
\end{equation}

Note that identities~\eqref{equ.s6rsk7b} and \eqref{equ.zxs4hwa} are special cases of identity~\eqref{equ.t70h9q0} at $j=1$ and $j=2$, respectively.

\bigskip

At $k=j\ne 0$ in identity~\eqref{equ.bxwz09x} we have

\begin{equation}\label{equ.izw0y7l}
\sum_{r=1}^\infty\tan^{-1}\left\{\frac{F_{2j}^2L_{4jr+2j}}{F_{4jr+2j}^2}\right\}=\tan^{-1}\left(\frac{1}{L_{2j}}\right)\,,
\end{equation}

yielding at $j=1$, the special value

\begin{equation}\label{equ.hqdqegn}
\boxed{\sum_{r=1}^\infty\tan^{-1}\left\{\frac{L_{4r+2}}{F_{4r+2}^2}\right\}=\tan^{-1}\left(\frac{1}{3}\right)}\,.
\end{equation}

\bigskip

Finally, taking limit of identity~\eqref{equ.bxwz09x} as $j\to\infty$, we obtain

\begin{equation}\label{equ.yb6rk2c}
\lim_{j\to\infty}\;\sum_{r=1}^\infty\tan^{-1}\left\{\frac{F_{2j}^2L_{4jr+2k}}{F_{4jr+2k}^2}\right\}=\tan^{-1}\left(\frac{1}{\phi^{2k}}\right)\,.
\end{equation}

\subsubsection{ $5\lambda^2=L_{4j}$, $p=0$ and $k=j$ in identity~\eqref{equ.o14ht74}}

Another interesting particular case of identity~\eqref{equ.o14ht74} is obtained by setting \mbox{$5\lambda^2=L_{4j}$}, $p=0$ and $k=j$ to obtain

\begin{equation}\label{equ.qzbk8sv}
\sum_{r = 1}^\infty  {\tan ^{ - 1} \left\{ {\frac{{F_{2j} \sqrt {5L_{4j} } L_{4jr - 2j} }}{{L_{2(4jr - 2j)} }}} \right\}}  = \frac{\pi}{2}\,,
\end{equation}

which at $j=1$ gives the special value

\begin{equation}
\boxed{\sum_{r = 1}^\infty  {\tan ^{ - 1} \left\{ {\frac{\sqrt{35}\, L_{4r - 2} }{{L_{2(4r - 2)} }}} \right\}}  = \frac{\pi}{2}}\,.
\end{equation}

\subsubsection{ $5\lambda^2=L_{4j}$, $p=0$ and $k=2j$ in identity~\eqref{equ.o14ht74}}

In this case Theorem~\ref{thm.xrc79et} reduces to

\begin{equation}
\sum_{r = 1}^\infty  {\tan ^{ - 1} \left\{ {\frac{{F_{2j} \sqrt {5L_{4j} } L_{4jr} }}{{L_{8jr} }}} \right\}}  = \tan ^{ - 1} \left( {\frac{1}{{\sqrt 5 }}\frac{{\sqrt {L_{4j} } }}{{F_{2j} }}} \right)\,.
\end{equation}

At $j=1$, we have the special value

\begin{equation}
\boxed{\sum_{r = 1}^\infty  {\tan ^{ - 1} \left\{ {\frac{\sqrt{35}\, L_{4r} }{{L_{8r} }}} \right\}}  = \sqrt{\frac{7}{5}}}\,.
\end{equation}

\subsubsection{$\lambda=L_{2j}/\sqrt 5$ and $k=j$ in identity~\eqref{equ.o14ht74}}

Setting $\lambda=L_{2j}/\sqrt 5$ and $k=j$ in identity~\eqref{equ.o14ht74} we have

\begin{equation}\label{equ.eyh9qun}
\sum_{r = p}^\infty  {\tan ^{ - 1} \left\{ {\frac{{\sqrt 5 F_{4j} }}{{L_{4jr + 2j}}}} \right\}}  = \tan ^{ - 1} \left( {\frac{{ L_{2j} }}{{F_{4jp}\sqrt 5 }}} \right)\,,
\end{equation}

which at $p=1$ gives

\begin{equation}\label{equ.m0qe2qb}
\sum_{r = 1}^\infty  {\tan ^{ - 1} \left\{ {\frac{{\sqrt 5 F_{4j}}}{{L_{4jr + 2j}}}} \right\}}  = \tan ^{ - 1} \left( {\frac{{1}}{{F_{2j}\sqrt 5 }}} \right)
\end{equation}

and at $p=0$ yields
\begin{equation}\label{equ.nwoiutl}
\sum_{r = 1}^\infty  {\tan ^{ - 1} \left\{ {\frac{{\sqrt 5 F_{4j}}}{{L_{4jr - 2j}}}} \right\}}  = \frac{\pi}{2}\,.
\end{equation}

\subsubsection{$\lambda=L_{2j}/\sqrt 5$, $p=0$ and $k=2j\ne 0$ in identity~\eqref{equ.o14ht74}}

The above choice yields
\begin{equation}\label{equ.lzwycuj}
\sum_{r = 1}^\infty  {\tan ^{ - 1} \left\{ {\frac{{\sqrt 5 F_{4j} }}{{L_{4jr}}}} \right\}}  = \tan ^{ - 1} \left( {\frac{{ L_{2j} }}{{F_{2j}\sqrt 5 }}} \right)\,.
\end{equation}

\subsection{Results from Theorem~\ref{thm.kc5pf5v}}

\subsubsection{$\lambda=F_{2j-1}$ and $p=1$ in identity~\eqref{equ.wz6niw6}}

The above choice gives

\begin{equation}\label{equ.mj4nx98}
\sum_{r = 1}^\infty  {\tan ^{ - 1} \left\{ {\frac{{F_{2(2j - 1)} }}{{F_{4jr - 2r + 2k - 1} }}} \right\}}  = \tan ^{ - 1} \left( {\frac{{F_{2j - 1} }}{{F_{2j + 2k - 2} }}} \right)\,.
\end{equation}

At $k=j$ in identity~\eqref{equ.mj4nx98} we have the interesting formula 

\begin{equation}\label{equ.tbzc0kv}
\boxed{\sum_{r = 1}^\infty  {\tan ^{ - 1} \left\{ {\frac{{F_{2(2j - 1)} }}{{F_{4jr - 2r + 2j - 1} }}} \right\}}  = \tan ^{ - 1} \left( {\frac{{1}}{{L_{2j - 1} }}} \right)}\,.
\end{equation}

Note that identity~\eqref{equ.tbzc0kv}, at $j=1$, includes Lehmer's result (cited in \cite{hoggattv, melham95}) as a particular case.

\bigskip

Setting $j=1$ in identity~\eqref{equ.mj4nx98} we obtain

\begin{equation}\label{equ.dsi0drc}
\boxed{\sum_{r = 1}^\infty  {\tan ^{ - 1} \left\{ {\frac{{1}}{{F_{2r + 2k - 1} }}} \right\}}  = \tan ^{ - 1} \left( {\frac{1}{{F_{2k} }}} \right)}\,.
\end{equation}

Note again that identity~\eqref{equ.dsi0drc} subsumes Lehmer's formula and the result of Melham ($p=1$ in identity(3.5) of \cite{melham95}), at $k=1$ and at $k=0$ respectively.

\bigskip

Finally, taking limit $j\to\infty$ in identity~\eqref{equ.mj4nx98}, we obtain
\begin{equation}\label{equ.cr0s48s}
\lim_{j\to\infty}\,\sum_{r = 1}^\infty  {\tan ^{ - 1} \left\{ {\frac{{F_{2(2j - 1)} }}{{F_{4jr - 2r + 2k - 1} }}} \right\}}  = \tan ^{ - 1} \left( {\frac{{1 }}{{\phi^{2k-1} }}} \right)\,.
\end{equation}

\subsubsection{$\lambda=L_{2j-1}/\sqrt 5$ and $k=j$ in identity~\eqref{equ.wz6niw6}}

The above choice gives

\begin{equation}\label{equ.a1f1zqv}
\sum_{r = p}^\infty  {\tan ^{ - 1} \left\{ {\frac{{\sqrt 5 L_{2j - 1}^2 F_{4jr - 2r + 2j - 1} }}{{L_{4jr - 2r + 2j - 1}^2 }}} \right\}}  = \tan ^{ - 1} \left( {\frac{1}{{\sqrt 5 }}\frac{{L_{2j - 1} }}{{F_{4jp - 2p} }}} \right)\,.
\end{equation}

Setting $p=1$ in identity~\eqref{equ.a1f1zqv}, we find

\begin{equation}\label{equ.xr9ngzb}
\sum_{r = 1}^\infty  {\tan ^{ - 1} \left\{ {\frac{{\sqrt 5 L_{2j - 1}^2 F_{4jr - 2r + 2j - 1} }}{{L_{4jr - 2r + 2j - 1}^2 }}} \right\}}  = \tan ^{ - 1} \left( {\frac{1}{{\sqrt 5 }}\frac{{1}}{{F_{2j - 1} }}} \right)\,,
\end{equation}

while choosing $j=1$ leads to

\begin{equation}\label{equ.awptwpq}
\sum_{r = p}^\infty  {\tan ^{ - 1} \left\{ {\frac{{\sqrt 5 F_{2r + 1} }}{{L_{2r + 1}^2 }}} \right\}}  = \tan ^{ - 1} \left( {\frac{1}{{\sqrt 5 }}\frac{{1}}{{F_{2p} }}} \right)\,,
\end{equation}

which at $p=0$ gives the special value

\begin{equation}\label{equ.v6sfa8i}
\boxed{\sum_{r = 1}^\infty  {\tan ^{ - 1} \left\{ {\frac{{\sqrt 5 F_{2r - 1} }}{{L_{2r - 1}^2 }}} \right\}}  = \frac{\pi}{2}}\,.
\end{equation}

\subsubsection{$5\lambda^2=L_{4j-2}$ and $k=j$ in identity~\eqref{equ.wz6niw6}}

The above substitutions give

\begin{equation}\label{equ.xdqih4r}
\sum_{r = p}^\infty  {\tan ^{ - 1} \left\{ {\frac{{\sqrt{5L_{4j - 2}} L_{2j - 1} F_{4jr - 2r + 2j - 1} }}{L_{2(4jr - 2r + 2j - 1)} }} \right\}}  = \tan ^{ - 1} \left( {\frac{\sqrt{5L_{4j - 2}} }{{5F_{4jp - 2p} }}} \right)\,.
\end{equation}

At $p=0$ in identity~\eqref{equ.xdqih4r} we have, for positive integers $j$,

\begin{equation}\label{equ.iwv0r4l}
\sum_{r = 1}^\infty  {\tan ^{ - 1} \left\{ {\frac{{\sqrt{5L_{4j - 2}} L_{2j - 1} F_{4jr - 2r - 2j + 1} }}{L_{2(4jr - 2r - 2j + 1)} }} \right\}}  = \frac{\pi}{2}\,,
\end{equation}

giving, at $j=1$, the special value

\begin{equation}\label{equ.v7chd82}
\boxed{\sum_{r = 1}^\infty  {\tan ^{ - 1} \left\{ {\frac{{\sqrt {15} F_{2r - 1} }}{{L_{2({2r - 1})} }}} \right\}}  = \frac{\pi}{2}}\,.
\end{equation}

At $p=2$ in identity~\eqref{equ.xdqih4r} we have, for positive integers $j$,

\begin{equation}\label{equ.mmnmget}
\sum_{r = 1}^\infty  {\tan ^{ - 1} \left\{ {\frac{{\sqrt{5L_{4j - 2}} L_{2j - 1} F_{4jr-2r+6j-3} }}{L_{2(4jr-2r+6j-3)} }} \right\}}  = \tan ^{ - 1} \left( {\frac{1}{{\sqrt{5F_{4j - 2}F_{8j - 4}} }}} \right)\,,
\end{equation}

which gives, at $j=1$, the special value

\begin{equation}
\boxed{\sum_{r = 1}^\infty  {\tan ^{ - 1} \left\{ {\frac{{\sqrt {15} F_{2r + 3} }}{{L_{2({2r + 3})} }}} \right\}}  = \tan ^{ - 1} \left( \frac{1}{\sqrt{15}} \,\right)}\,.
\end{equation}

\subsection{Results from Theorem~\ref{thm.xi2xcnh}}

\subsubsection{$\lambda=\sqrt {L_{4j}}$, $k=0$ and $p=1$ in identity~\eqref{equ.ywa2m9d}}

The above choice gives

\begin{equation}\label{equ.qk4nj3o}
\sum_{r=1}^\infty\tan^{-1}\left(\frac{5\sqrt {L_{4j}} F_{2j}F_{4jr-1}}{L_{8jr-2}}\right)=\tan^{-1}\left(\frac{\sqrt {L_{4j}}}{L_{2j-1}}\right)\,,
\end{equation}

which, at $j=1$, gives

\begin{equation}\label{equ.s9mmd68}
\boxed{\sum_{r=1}^\infty\tan^{-1}\left(\frac{5\sqrt {7} F_{4r-1}}{L_{2(4r-1)}}\right)=\tan^{-1}\sqrt 7}\,.
\end{equation}

\subsubsection{$\lambda=L_{2j}$ and $p=1$ in identity~\eqref{equ.ywa2m9d}}

Setting $\lambda=L_{2j}$ and $p=1$ in identity~\eqref{equ.ywa2m9d} gives

\begin{equation}\label{equ.ksxz0op}
\sum_{r=1}^\infty \tan^{-1}\left\{\frac{F_{4j}}{F_{4jr+2k-1}}\right\}=\tan^{-1}\left(\frac{L_{2j}}{L_{2j+2k-1}}\right)\,.
\end{equation} 

\bigskip

Taking limit as $j\to\infty$ in identity~\eqref{equ.ksxz0op} gives
\begin{equation}\label{equ.b2gtlh1}
\lim_{j\to\infty}\sum_{r=1}^\infty \tan^{-1}\left\{\frac{F_{4j}}{F_{4jr+2k-1}}\right\}=\tan^{-1}\left(\frac{1}{\phi^{2k-1}}\right)\,.
\end{equation}

\subsubsection{$\lambda=\sqrt 5F_{2j}$, $p=1$ and $k=0$ in identity~\eqref{equ.ywa2m9d}}

Setting $\lambda= {\sqrt 5F_{2j}}$, $p=1$ and $k=0$ in identity~\eqref{equ.ywa2m9d} we obtain

\begin{equation}\label{equ.r677iuu}
\sum_{r=1}^\infty\tan^{-1}\left(\frac{5\sqrt 5 {F_{2j}^2}F_{4jr-1}}{L_{4jr-1}^2}\right)=\tan^{-1}\left(\frac{{\sqrt 5F_{2j}}}{L_{2j-1}}\right)\,,
\end{equation}

which gives the special value

\begin{equation}\label{equ.qplrlge}
\sum_{r=1}^\infty\tan^{-1}\left(\frac{5\sqrt 5 F_{4r-1}}{L_{4r-1}^2}\right)=\tan^{-1}\sqrt 5\,,
\end{equation}

at $j=1$.

\subsection{Results from Theorem~\ref{thm.vqhpot1}}

\subsubsection{$\lambda=\sqrt {L_{4j-2}}$ and $j=0=k$ in identity~\eqref{equ.ocyrpi4}}

With the above choice we obtain

\begin{equation}\label{equ.dshfqze}
\sum_{r = p}^\infty  {\tan ^{ - 1} \left\{ {\frac{{\sqrt 3 L_{2r} }}{{L_{4r} }}} \right\}}  = \tan ^{ - 1} \left( {\frac{{\sqrt 3 }}{{L_{2p - 1} }}} \right)\,,
\end{equation}

which gives rise, at $p=1$, to the special value

\begin{equation}
\boxed{\sum_{r = 1}^\infty  {\tan ^{ - 1} \left\{ {\frac{{\sqrt 3 L_{2r} }}{{L_{4r} }}} \right\}}  = \frac{\pi }{3}}\,.
\end{equation}

\subsubsection{$\lambda= {L_{2j-1}}$ and $p=1$ in identity~\eqref{equ.ocyrpi4}}

With the above choice we have

\begin{equation}\label{equ.ke7jl7n}
\sum_{r = 1}^\infty  {\tan ^{ - 1} \left\{ {\frac{{L_{2j - 1}^2 }}{5}\frac{{L_{4jr - 2r + 2k} }}{{F_{4jr - 2r + 2k}^2 }}} \right\}}  = \tan ^{ - 1} \left( {\frac{{L_{2j - 1} }}{{L_{2j + 2k - 1} }}} \right)\,.
\end{equation}

$k=0$ in identity~\eqref{equ.ke7jl7n} gives 

\begin{equation}\label{equ.awl45g6}
\sum_{r = 1}^\infty  {\tan ^{ - 1} \left\{ {\frac{{L_{2j - 1}^2 }}{5}\frac{{L_{4jr - 2r} }}{{F_{4jr - 2r}^2 }}} \right\}}  =   {\frac{\pi}{4}}\,,
\end{equation}

which at $j=1$ gives the special value
\begin{equation}\label{equ.ny4a2vn}
\boxed{\sum_{r = 1}^\infty  {\tan ^{ - 1} \left\{ {\frac{1}{5}\frac{{L_{2r} }}{{F_{2r}^2 }}} \right\}}  = \frac{\pi }{4}}\,.
\end{equation}

$j=1$ in identity~\eqref{equ.ke7jl7n} leads to
\begin{equation}
\sum_{r = 1}^\infty  {\tan ^{ - 1} \left\{ {\frac{1}{5}\frac{{L_{2r+2k} }}{{F_{2r+2k}^2 }}} \right\}}  = \tan^{-1}\left(\frac{1}{L_{2k+1}}\right)\,,
\end{equation}

which gives the special value

\begin{equation}\label{equ.le9h1ec}
\boxed{\sum_{r = 1}^\infty  {\tan ^{ - 1} \left\{ {\frac{1}{5}\frac{{L_{2r+2} }}{{F_{2r+2}^2 }}} \right\}}  = \tan^{-1}\left(\frac{1}{4}\right)}\,,
\end{equation}

at $k=1$.

\bigskip

Taking limit $j\to\infty$ in identity~\eqref{equ.ke7jl7n}, we obtain

\begin{equation}\label{equ.lf3mqef}
\lim_{j\to\infty}\,\sum_{r = 1}^\infty  {\tan ^{ - 1} \left\{ {\frac{{L_{2j - 1}^2 }}{5}\frac{{L_{4jr - 2r + 2k} }}{{F_{4jr - 2r + 2k}^2 }}} \right\}}  = \tan ^{ - 1} \left( {\frac{1}{\phi^{2k}}} \right)\,.
\end{equation}

\subsubsection{$\lambda= {L_{2j-1}}$ and $j=0=k$ in identity~\eqref{equ.ocyrpi4}}

This choice gives

\begin{equation}\label{equ.v7p9pwu}
\sum_{r = p}^\infty  {\tan ^{ - 1} \left\{ {\frac{1}{5}\frac{{L_{2r} }}{{F_{2r}^2 }}} \right\}}  = \tan ^{ - 1} \left( {\frac{1}{{L_{2p - 1} }}} \right)\,,
\end{equation} 

Note that identities~\eqref{equ.ny4a2vn} and \eqref{equ.le9h1ec} are special cases of \eqref{equ.v7p9pwu} at $p=1$ and at $p=2$.

\subsubsection{$\lambda=\sqrt 5F_{2j-1}$ and $j=0=k$ in identity~\eqref{equ.ocyrpi4}}

The above choice gives

\begin{equation}\label{equ.dm04elr}
\sum_{r = p}^\infty  {\tan ^{ - 1} \left\{ {\frac{{\sqrt 5 }}{{L_{2r} }}} \right\}}  = \tan ^{ - 1} \left( {\frac{{\sqrt 5 }}{{L_{2p - 1} }}} \right)\,,
\end{equation}

which at $p=1$ gives the special value

\begin{equation}
\boxed{\sum_{r = 1}^\infty  {\tan ^{ - 1} \left\{ {\frac{{\sqrt 5 }}{{L_{2r} }}} \right\}}  = \tan ^{ - 1} \sqrt 5}\,. 
\end{equation}

\section{Conclusion}

Using a fairly straightforward technique, we have derived numerous infinite arctangent summation formulas involving Fibonacci and Lucas numbers. While most of the results obtained are new, a couple of `celebrated' results appear as particular cases of more general formulas derived in this paper.

\end{document}